\input amstex
\input prepictex
\input pictex
\input postpictex
\magnification=\magstep1
\documentstyle{amsppt}
\TagsOnRight
\hsize=5.1in                                                  
\vsize=7.8in
\define\R{{\bold R}}
\define\cl{\operatorname {cl}}
\def\k{\bold k}

\vskip 2cm
\topmatter

\title Topology of closed 1-forms and their critical points \endtitle

\rightheadtext{}
\leftheadtext{}
\author  Michael Farber \endauthor
\address
School of Mathematical Sciences,
Tel-Aviv University,
Ramat-Aviv 69978, Israel
\endaddress
\email farber\@math.tau.ac.il
\endemail
\thanks{The research was supported by a grant from the 
Israel Academy of Sciences and Humanities and by
the Herman Minkowski Center for Geometry.} 
\endthanks
\abstract{In this paper we suggest an analog of the Lusternik - Schnirelman 
theory for closed 1-forms. Namely, we use cup-products and higher Massey products
to find topological lower bounds on the minimal number of geometrically distinct critical points
of any closed 1-form in a given cohomology class.}
\endabstract
\endtopmatter
%-------------------------------------------------------------------------
%----------------------------------------------------------------------

\define\C{{\bold C}}

\define\Z{{\bold Z}}
\define\Q{{\bold Q}}
\define\V{{\Cal V}}

\define\Hom{\operatorname{Hom}}

\define\Tor{\operatorname{Tor}}

\define\rk{\operatorname{rk}}

\redefine\c{{\frak c}}
%\define\cl{{\frak Cl}}

\define\im{\operatorname{im}}

%\define\Ker{\operatorname{Ker}}

\define\rank{\operatorname{rank}}
%\define\ln{\operatorname{ln}}

\define\E{\Cal E}

\def\<{\langle}
\def\>{\rangle}

\define\pd#1#2{\dfrac{\partial#1}{\partial#2}}

\def\part{\partial}

\def\lam{\lambda}

\NoBlackBoxes

\def\cat{\operatorname {cat}}

\def\pn{\partial_+N}

\def\fp{\Bbb F_p}

\documentstyle{amsppt}

\nopagenumbers

\heading{\bf \S 1. Introduction} \endheading

Let $X$ be a closed manifold and let $\xi\in H^1(X;\R)$ be a nonzero cohomology class. 
The well-known Novikov inequalities \cite{N1, N2} estimate the numbers of critical points of
different indices of any closed 1-form $\omega$ on $X$ lying in the class $\xi$, assuming that all the
singular points are non-degenerate in the sense of Morse. 
Novikov type inequalities were generalized in \cite{BF1} for closed 1-forms 
with more general singularities (non-degenerate in the sense of Bott). 
In \cite{BF2} an equivariant generalization of the Novikov inequalities was developed. 

Novikov inequalities have found interesting and important applications in symplectic topology, especially in the study of 
symplectic fixed points (Arnold's conjecture). 
Here we should mention the work of J.-C. Sikorav \cite{S}, Hofer - Salamon \cite{HS}, Van - Ono \cite{VO}, and most
recently the preprint of Eliashberg and Gromov \cite{EG}. 

In this paper we describe new results, 
which give topological restrictions on the number of geometrically
distinct critical points of closed 1-form in a given cohomology class. We  impose no assumptions
on the nature of the critical points. Therefore, the results of this paper have the same 
relation to the classical Lusternik - Schnirelman - Frolov - Elsgoltz theory, as Novikov's theory has
to the classical Morse theory.

The main theorem of the paper states that any closed 1-form $\omega$ on a closed $n$-dimensional manifold must have at least $\cl_\k(\xi)-1$ geometrically distinct critical points, where
$\xi\in H^1(M;\R)$ is the cohomology class of $\omega$. Here $0\le \cl_\k(\xi)\le n$ is the number,
which we define and study in this paper; we 
call it {\it the cup-length associated with the class $\xi$}. It is defined using the 
cohomological cup-products in flat bundles, which are somehow related to $\xi$ 
(are {\it $\xi$-generic}). We prove also a theorem, which produces easily computable
lower bounds on the
number of critical points of closed 1-forms using higher
Massey products (instead of $\xi$-generic bundles).

We show by example that our estimates prove that in some cases
any closed 1-form on $M$ has at least $n-1$ critical points, where $n=\dim M$. This type of 
estimate
cannot be improved since, as we show in this paper using the methods of Takens \cite{T}, 
on any closed 
$n$-dimensional manifold, in any integral cohomology class $\xi$,
there always exists a closed 1-form with $\le n-1$ critical points.

The main technical tool of our proof is {\it the deformation complex} which we describe in 6.3.
It provides a way of dealing simultaneously and polynomially with all flat bundles of the form $a^\xi\otimes F$, where $a\in \k$, and also with their "limit as $a\to 0$".

Different estimates on the number of critical points of closed 1-forms were recently suggested
in \cite{F2, F3}, where we used flat line bundles which are described by complex numbers, which are not Dirichlet units. 
The approach of the present paper has some important advantages:
simplicity, a larger variety of flat bundles which can be used to produce the estimates, the possibility of using
different fields (for example, fields having positive characteristic). There are examples when the approach of this paper
gives stronger estimates than the approach of \cite{F2, F3}, although in some cases the situation
is the opposite.

The approach developed in the present
paper, because of its simplicity and flexibility, 
could have immediate infinite dimensional generalizations (for closed 1-forms on Banach manifolds 
with Palais - Smale type conditions). Hence, the methods of this paper may be applied
to study of the problem 
of estimating the number of closed trajectories of Hamiltonian systems, which was the
main motivation of S. P. Novikov in \cite{N1, N2}, while developing his theory. 
One may also hope to develop a Floer theory version
of the main theorem of the present paper, which would produce applications in the theory
of symplectic fixed points. 

I am grateful to Shmuel Weinberger for his useful comments.

\heading{\bf \S 2. The cup-length associated with a cohomology class}\endheading
 
\subheading{2.1. Notation} 
Let $\k$ be a fixed algebraically closed field. The most important cases, which the reader should
keep in mind are $\k=\C$ or $\k$ being the algebraic closure of a finite field $\fp$.

We will consider flat $\k$-vector bundles $E$ over a compact polyhedron $X$. 
We will understand such bundles as
locally trivial sheaves of $\k$-vector spaces. The cohomology $H^q(X;E)$ will be understood as the
sheaf cohomology. 

A flat vector bundle is
determined by its monodromy -- linear representation of the fundamental group $\pi_1(X,x_0)$ 
on the fiber $E_0$ over the base point $x_0$, which is given by the parallel transport along loops.
For example, a flat $\k$-line bundle is determined by a homomorphism $H_1(X)\to \k^\ast$; here $\k^\ast$ is considered 
as a multiplicative abelian group.

Given a real cohomology class $\xi\in H^1(X;\R)$, it can be viewed as a homomorphism
$\xi: H_1(X)\to \R$, and we will denote the kernel by $\ker (\xi)$. Let $\V_\xi$ be
the variety of all $\k$-line bundles over $X$, 
which have trivial monodromy along the curves in $\ker(\xi)$.
If $r$ denotes the rank 
$r \, =\, \rank (H_1(X)/\ker (\xi))$, then $\V_\xi$ can be identified with $(\k^\ast)^r$.

\subheading{2.2. Definition} {\it
A flat bundle $F$ will be called $\xi$-generic if there is no $E\in \V_\xi$, so that for some $q$,
$\dim H^q(X;E\otimes F) < \dim H^q(X;F)$.}

Note that this property depends only on $\ker(\xi)$.

It follows that if we have two classes $\xi, \eta\in H^1(X;\R)$ and $\ker(\xi) \subset \ker(\eta)$
then $\V_\eta\subset \V_\xi$ and {\it any flat bundle, which is $\xi$-generic is also $\eta$-generic.}

The examples of $\xi$-generic bundles can be constructed as follows. Suppose for simplicity
that the class 
$\xi$ is integral, $\xi\in H^1(X;\Z)$. Any $a\in \k^\ast$ determines a flat $\k$-line bundle over $X$,
so that the monodromy along any loop $\gamma\in \pi_1(X)$ is given by 
$a^{\<\xi,\gamma\>}\in \k^\ast$. We will denote this bundle by $a^\xi$. It is clear that the bundle $a^\xi$ is $\xi$-generic for
almost all $a$, i.e. for all, except finitely many. For instance, in the case $\k=\C$ the bundle $a^\xi$
is $\xi$-generic for all {\it transcendental} $a\in \C$, cf. \cite{F3}.

A picture illustrating the notion of a $\xi$-generic bundle, gives Proposition 2.3 below. 
This statement describes the behavior of the function 
$E\mapsto \dim_\k H^i(X;E\otimes F)$, when $E$ runs over flat $\k$-line bundles $E\in \V_\xi$
and $F\to X$ is an arbitrary fixed flat $\k$-vector bundle. 
\proclaim{2.3. Proposition} Let $X$ be a compact polyhedron, $F\to X$ be a flat $\k$-vector
bundle, and
$\xi\in H^1(X;\R)$ be a real cohomology class. 
For a fixed $q$, the function
$$E\mapsto \dim_\k H^q(X;E\otimes F),\qquad E\in \V_\xi,\tag2-1$$
has the following behavior:

(a)  it is constant for all $E$ lying outside an algebraic subvariety
$\Sigma_q(F) \subset \V_\xi$ (the jump subvariety);

(b) for $E\in \Sigma_q(F)$, 
the dimension of the cohomology $\dim_\k H^q(X;E\otimes F)$ is greater than the dimension of 
$\dim_\k H^q(X;E\otimes F)$ for $E\in (\V_\xi-\Sigma_q(F))$;

(c) in the case $\k =\C$ and $F$ being the trivial line bundle,
the common value of $\dim_\k H^q(X;E)$ for $E\in (\V_\xi -\Sigma_q(F))$ equals the Novikov 
number $b_q(\xi)$.
\endproclaim

We refer to \cite{N1-N3} for the definition of the Novikov numbers associated with the cohomology class
$\xi\in H^1(X;\R)$.
We will also use the definition suggested in \cite{F} (using rings of rational functions instead of
formal power series), which (as shown in \cite{F}) is
equivalent to Novikov's original definition.

Proposition 2.3 is well known.

\subheading{2.4. Definition of the cup-length $\cl_\k (\xi)$} {\it
Given a cohomology class $\xi\in H^1(X;\R)$ we will define
the cup-length associated with $\xi$ (denoted $\cl_{\k}(\xi)$) as the largest integer 
$m= 0, 1,2, ...$, so that there exist flat $\k$-vector bundles $E_1, E_2, \dots, E_m$
and cohomology classes 
$$v_1\in H^{d_1}(X;E_1), v_2\in H^{d_2}(X;E_2), \dots,
v_m\in H^{d_m}(X;E_m),\tag2-2$$
such that:
\newline(i) $E_1$ and $E_2$ are $\xi$-generic;\newline
(ii) $d_i>0$ for $i=3, \dots, m$;\newline
(iii)  the cup-product 
$$v_1\cup v_2\cup v_3\cup\dots\cup v_m \, \in H^d(X;E)\tag2-3$$
is nontrivial.} Here $d=\sum d_i$ and $E=E_1\otimes E_2\otimes \dots\otimes E_m$.

We will not focus our attention on the cases when $\cl_\k (\xi)\le 1$ since then the main Theorem 3.1
will be vacuous.

In general, the numbers $d_1$ and $d_2$ (which appear in the above definition)
are allowed to be zero. For example, for $\xi=0$ we may always take $E_1=E_2=\k$ (the trivial line bundle) and $v_1=v_2=1\in H^0(X;\k)$. This shows that $\cl_\k(0)$ {\it is the usual cup-length of $X$ plus 2.}

Note that $d_1$ and $d_2$ are automatically positive, assuming that the class $\xi\ne 0$
is nontrivial and $X$ is connected. Hence we obtain an inequality
$$\cl_{\k}(\xi)\le \dim(X),\qquad \xi\ne 0.\tag2-4$$

Here is another observation.
Suppose that $X$ is a closed orientable manifold.
{\it Then $\cl_{\k}(\xi)\ge 2$ if and only if some cohomology
$H^\ast(X;F)$ is non-trivial for a $\xi$-generic flat $\k$-vector bundle 
$F$.} Indeed, 
if $H^i(X;F)$ is nontrivial then by the Poincar\'e duality there exists a nontrivial cup-product
$H^i(X;F)\otimes H^{n-i}(X;F^{\ast})\to \k,$
where $n=\dim X$ and $F^{\ast}$ is the dual flat $\k$-vector bundle. 
It is clear that we may choose a $\xi$-generic $F$ so that $F^{\ast}$ is also $\xi$-generic.

In the case $\k=\C$, a slightly weaker statement says that
$\cl_\k(\xi)\ge 2$ if at least one of the Novikov numbers $b_i(\xi)$ is positive. 
This follows from Proposition 2.3.(c).

\heading{\bf \S 3. The main result}\endheading

\proclaim{3.1. Main Theorem} Let $X$ be a closed manifold and
let $\omega$ be a closed 1-form on $X$ lying in a cohomology class
$\xi\in H^1(X;\R)$. Let $S(\omega)$ denote the set of critical points of $\omega$, i.e. the set of points $p\in X$
such that $\omega_p=0$. Then the  Lusternik - Schnirelman category of $S(\omega)$ satisfies
$$\cat(S(\omega)) \ge \cl_\k(\xi)-1.\tag3-1$$
In particular, if the set of critical points $S(\omega)$ is finite 
then for the total number $\#S(\omega)$ of critical points,
$$\#S(\omega)\, \ge \, \cl_\k(\xi)-1.\tag3-2$$
\endproclaim

Here $\cat(S)$ denotes the classical Lusternik - Schnirelman category of $S=S(\omega)$, i.e.
the least number $r$, so that $S$ can be covered by $r$ closed subsets 
$A_1\cup\dots \cup A_r$ such that each inclusion $A_j\to S$ is null-homotopic.

Since $\cl_\k(\xi)\ge n$, where $n=\dim X$, the highest possible lower bound given by Theorem 3.1
is $\# S(\omega) \ge n-1$. In the next section we will see examples where any closed 1-form on
a given closed manifold of dimension $n$ has at least $n-1$ critical points. This estimate is the highest
possible, as shown by the following result:

\proclaim{3.2. Theorem} Let $X$ be a closed connected $n$-dimensional manifold, and let 
$\xi\in H^1(X;\Z)$ be a nonzero integral 
cohomology class. Then there exists a closed 1-form $\omega$
on $X$, realizing $\xi$, having at most $n-1$ critical points. \endproclaim

The proof is given in section \S 8, using the method of Takens \cite{T}.

\subheading{3.3. Application: zeros of symplectic vector fields} As an immediate application 
of the main Theorem 3.1 we will mention here a topological estimate of the number of geometrically
distinct zeros which must have a symplectic vector field.

Let $X$ be a symplectic manifold
with the symplectic form $\nu$. Recall that a smooth vector field $V$ on $X$ is called {\it symplectic}
if the Lie derivative of $\nu$ with respect to $V$ vanishes, i.e. $\Cal L_V(\nu)=0$. 

\proclaim{Corollary} Let $V$ be a symplectic vector field, and let $S=S(V)$ denote the set of zeros
of $V$. Then the Lusternik - Schnirelman category of $S$ satisfies
$$\cat(S)\ge \cl_\k(\xi_V)-1,\tag3-3$$
where $\xi_V\in H^1(X;\R)$ denotes the cohomology class of the closed 1-form $\omega=i_V(\nu)$.
In particular, if $S$ is finite, then for the number of distinct zeroes of $V$ we have
$$\#(S)\ge \cl_\k(\xi_V)-1.\tag3-4$$
\endproclaim

\heading{\bf \S 4. Examples}\endheading

\subheading{4.1} Here we will construct in each dimension $n\ge 2$ a closed 
$n$-dimensional manifold
$X$ with non-trivial first Betti number so that 
{\it any closed 1-form on $X$ has at least $n-1$ critical points.}

Let $X$ be the connected sum of the real projective space $\R\bold P^n$
and $S^1\times S^{n-1}$. In other words,
$X$ is obtained from $\R\bold P^n$ by adding a handle of index one. 
Real cohomology $H^1(X;\R)$ is one dimensional, generated by an integral class $\xi$
which is trivial on the projective space and goes along the handle once. 
We will show using Theorem 3.1 that any closed 1-form in class $\xi$ has at least $n-1$ critical 
points. This would imply our above 
statement since for any closed 1-form $\omega$ on $X$ 
lying in a nontrivial cohomology class we will have $\lambda\omega$ representing $\xi$ for some $\lambda\in\R$.
Also, from our arguments below it will be
clear that the usual cup-length of $X$ equals $n$ and so any function on $X$
must have at least $n+1$ critical points, according to the Lusternik - Schnirelman theory.

In order to apply Theorem 3.1 we have first to choose a field $\k$. We take $\k$ to be the algebraic
closure of the field $\Z_2$. Flat $\k$-line bundles $E$ over $X$ belonging to $\V_\xi$ must be trivial
over the projective space and they are determined by their monodromy $a\in \k^\ast$ along the handle.
Such a line bundle is {\it generic} if and only if $a\ne 1$. 

Let $a\in k^\ast$ be an $n$-th root of 1, and let $E$ be a flat $\k$-line bundle over $X$ 
which is trivial over the projective space and has monodromy
$a$ along the handle. We will show that there is a class $u\in H^1(X;E)$ so that $0\ne u^n\in H^n(X;E^{\otimes n}) =H^n(X;\k)$.
This would give
$\cl_{\k}(\xi) \ge n,$
and so Theorem 3.1 implies our statement.

Let $\bar X$ denote the union of $\R\bold P^n$ and $S^1\times S^{n-1}$ intersecting in a small
$n$-dimensional disk; thus, $X$ is obtained from $\bar X$ by removing the interior of this disk.
It is clear that the flat bundle $E$ over $X$ extends to a flat bundle $\bar E$ over $\bar X$
so that $\bar E|_{X} \simeq E$. Since 
$$H^i(\bar X;\bar E) \simeq H^i(\R\bold P^n; \bar E|_{\R\bold P^n})\oplus 
H^i(S^1\times S^{n-1}; \bar E|_{S^1\times S^{n-1}})$$
for $i>0$, we find a class $\bar u\in H^1(\bar X;\bar E)$ so that
$${\bar u}^n\ne 0,\quad\text{and}\quad \bar u|_{S^1\times S^{n-1}}=0.$$
Let $u = \bar u|_X \in H^1(X;E)$. We claim that $u^n\ne 0$. Since the restriction map
$$H^n(\R\bold P^n;\k)\oplus H^n(S^1\times S^{n-1};\k)\simeq H^n(\bar X;\k)\to H^n(X;\k)$$
maps each of the summands isomorphically, we obtain that $u^n\ne 0$. \qed

\subheading{4.2} Consider now an example having a higher first Betti number and 
so admitting forms of higher
rank. Let 
$$X = \R\bold P^n\# (T^q\times S^{n-q}).$$
Arguments similar to those used above show that {\it any closed 1-form on $X$ has at least $n-1$
critical points.}

\subheading{4.3} Let $n$ be odd, $n=2\ell+1$, and $X$ be obtained by adding a handle of index 1 to 
a lens space $S^{2\ell+1}/\Z_p$, where $p$ is a prime number. Then for some $v\in H^1(X;\Z_p)$ and 
$u\in H^2(X;\Z_p)$, $v\cup u\cup u\cup \dots \cup u\ne 0$ ($\ell$ times
factor $u$). Hence, taking for $\k$ the algebraic closure of the finite Galois field $\Bbb F_p$, 
we obtain, using arguments similar to example 1, $\cl_\k(\xi) \ge \ell+1$ and so
$\#S(\omega) \ge \ell$
for any closed 1-form $\omega$ in class $\xi$, where $\xi \in H^1(X;\R)$ is a generator. 

Note the Theorem in \cite{F2, F3} gives no positive information in examples 4.1 - 4.3.

Note also that we 
had some other examples in \cite{F3}, which prove the existence of at least $n-1$ critical points
for closed 1-forms on $n$-dimensional manifolds. In example 3 in section 1.8 of \cite{F3} (the
3-manifold obtained by surgery on the knot $5_2$) the methods
of the present paper give no positive information, although the methods of \cite{F3} predict the existence
of critical points.

\heading{\bf \S 5. Main Theorem in terms of higher Massey products}\endheading

Our aim in this section 
is to find estimates on the cup-length $\cl_\k (\xi)$ 
using higher Massey products.

\subheading{5.1. Massey products} We will deal here with a special kind of higher Massey 
operations $d_r$ (where $r=1, 2, \dots$),
determined by a one-dimensional cohomology class $\xi\in H^1(X;\Z)$; these operations
were also described by S.P. Novikov \cite{N3}. The first operation $d_1: H^i(X;\k) \to 
H^{i+1}(X;\k) $ is the usual cup-product 
$$d_1(v)=(-1)^{i+1}v\cup\xi,\qquad\text{for}\quad v\in H^i(X;\k).\tag5-1$$
The higher Massey products $d_r: E_r^\ast \to E_r^{\ast+1}$ with $r\ge 2$, are defined as the differentials of 
a spectral sequence $(E_r^\ast, d_r)$, ($r\ge 1$)
with the initial term $E_1^\ast =H^\ast(X;\k)$ and the initial differential 
$d_1: E_1^\ast\to E_1^{\ast+1}$ given by (5-1). Each subsequent term $E_r^\ast$ is the cohomology of the previous
differential $d_{r-1}$: 
$$E_r^\ast=\ker(d_{r-1})/\im(d_{r-1}).\tag5-2$$
Traditionally, the following notation is used
$$d_r(v) = (-1)^{i+1}\cdot
\<v, \underset{\text{($r$ times)}}\to {\underbrace{\xi, \xi, \dots, \xi}}\>,\qquad
v\in H^i(X;\k).\tag5-3$$
We will describe this spectral sequence in full detail in section \S 9 below.

\subheading{Definition} A cohomology class $v\in H^i(X;\k)$ is said to be {\it $\xi$-surviver}
if $d_r(v)=0$ for all $r\ge 1$.

\proclaim{5.2. Theorem} Let $X$ be a closed manifold and $\xi\in H^1(X;\Z)$ be an integral 
cohomology class. Suppose that there exists a nontrivial cup-product
$$v_1\cup v_2\cup \dots\cup v_m\ne 0,\tag5-4$$
where the first two classes $v_1\in H^{d_1}(X;\k)$ and $v_2\in H^{d_2}(X;\k)$ are $\xi$-survivors 
and for $i=3, \dots, m$ the classes $v_i\in H^{d_i}(X;E_i)$ belong to the cohomology of some flat 
$\k$-vector bundles
$E_i$ over $X$ with $d_i>0$. Then
$$\cl_{k}(\xi) \ge m\tag5-5$$
and hence any closed 1-form $\omega$ in class $\xi$ has at least $m-1$ critical points.
\endproclaim

In case $\xi=0$ (when we study critical points of functions) the class $1\in H^0(X;\k)$ is a 
$\xi$-surviver. Hence in this case we may take $v_1=v_2=1$. This shows that in the case of functions
Theorem 5.2 is reduced to the usual Lusternik - Schnirelman inequality: we have $\cl_\k(\xi) = \cl(X)+2$
and so the number of critical points of a function on $X$ is $\ge  \cl_\k(\xi)- 1 \ge \cl(X)+1$ as the 
Lusternik - Schnirelman theory states.

The proof of Theorem 5.2 is given in \S 9.

\subheading{5.3. Detecting $\xi$-survivors} The following criterion allows us to show in some 
cases that a given cohomology class $v\in H^i(X;\k)$ is a $\xi$-surviver, where $\xi\in H^1(X;\Z)$.

{\it Suppose that we may realize $\xi$ by a smooth codimension one submanifold $V\subset X$,
having a trivial normal bundle, and we may realize the class $v$ by a simplicial cochain $c$,
such that the support of $c$ is disjoint from $V$. Then $v$ is a $\xi$-surviver.}

Indeed, this follows immediately from the definition of Massey operations $d_r$ (cf. (9-18) and (9-17))
and also from the precise formula for $\delta_1$ (given in \S 9.3).

This applies to the class $\xi$ itself and shows that it is a $\xi$-surviver, since we may realize it
by a cochain with support on a parallel copy of $V$. However,
this observation  is useless in producing nontrivial
products as in Theorem 5.2, because of the presence of another $\xi$-surviver, which kills the
whole product.

\subheading{5.4. Example} Let $X= \R\bold P^n\# (S^1\times S^{n-1})$ (the real 
projective space with a handle of index 1), and let $\xi\in H^1(X;\Z)$ be a class which 
restricts as the generator along the
handle and which is trivial on the projective space (as in 4.1). 
We may realize this class $\xi$ by a sphere $V=S^{n-1}$
cutting the handle. Let $\k$ be the algebraic closure of $\Z_2$. We have a class
$v\in H^1(X;\k)$, which restricts as an obvious generator of $H^1(\R\bold P^n;\k)$ and which is trivial 
on the handle. Applying the above criterion, we see that $v$ is a $\xi$-surviver 
(since we may realize it by a chain with support on the projective space, 
i.e. disjoint from $V$). As in 4.1 we obtain $v^n\ne 0$
and hence by Theorem 5.1, $\cl_\k(\xi)\ge n$. 

\subheading{5.5. Example} Let $\Sigma_g$ be a Riemann surface of genus $g>1$. Consider the 
classes $v_1, v_2, \xi\in H^1(\Sigma_g;\Z)$ which are Poincar\`e dual to the curves
shown in Figure 1.
\midinsert
$$\beginpicture
\setcoordinatesystem units <1.00000cm,1.00000cm>
\linethickness=1pt
\setlinear
%
% Fig CIRCULAR ARC object
%
\linethickness= 0.500pt
\circulararc 87.509 degrees from  3.302 21.145 center at  4.493 22.389
%
% Fig CIRCULAR ARC object
%
\linethickness= 0.500pt
\circulararc 50.506 degrees from  5.397 20.923 center at  4.443 18.987
%
% Fig CIRCULAR ARC object
%
\linethickness= 0.500pt
\circulararc 102.755 degrees from  7.588 21.336 center at  8.868 22.306
%
% Fig CIRCULAR ARC object
%
\linethickness= 0.500pt
\circulararc 46.809 degrees from  9.779 21.018 center at  8.858 18.891
%
% Fig CIRCULAR ARC object
%
\linethickness= 0.500pt
\circulararc 105.841 degrees from  4.381 22.796 center at  5.021 21.983
%
% Fig CIRCULAR ARC object
%
\linethickness= 0.500pt
\circulararc 102.072 degrees from  9.525 22.638 center at  9.782 21.638
%
% Fig ELLIPSE
%
\linethickness= 0.500pt
\ellipticalarc axes ratio  4.763:1.905  360 degrees 
	from 11.398 21.145 center at  6.636 21.145
%
% Fig ELLIPSE
%
\linethickness= 0.500pt
\ellipticalarc axes ratio  1.778:1.143  360 degrees 
	from  6.350 21.145 center at  4.572 21.145
\linethickness= 0.500pt
\setdots < 0.0953cm>
%
% Fig CONTROL PT SPLINE
%
% open spline
%
\plot	 4.413 22.860  4.572 22.701
 	 4.643 22.613
	 4.699 22.507
	 4.739 22.383
	 4.753 22.314
	 4.763 22.241
	 4.768 22.165
	 4.770 22.090
	 4.768 22.015
	 4.763 21.939
	 4.753 21.864
	 4.739 21.788
	 4.721 21.713
	 4.699 21.638
	 4.675 21.566
	 4.651 21.501
	 4.604 21.392
	 4.556 21.310
	 4.508 21.257
	 /
\plot  4.508 21.257  4.413 21.177 /
\linethickness= 0.500pt
%
% Fig CONTROL PT SPLINE
%
% open spline
%
\plot	 9.620 22.606  9.636 22.352
 	 9.638 22.288
	 9.635 22.221
	 9.628 22.153
	 9.616 22.082
	 9.600 22.010
	 9.580 21.935
	 9.555 21.859
	 9.525 21.780
	 9.492 21.704
	 9.455 21.632
	 9.414 21.565
	 9.370 21.503
	 9.272 21.394
	 9.160 21.304
	 /
\plot  9.160 21.304  8.922 21.145 /
%
% Fig TEXT object
%
\put{$v_2$} [lB] at  5.906 19.748
%
% Fig TEXT object
%
\put{$v_1$} [lB] at  4.477 23.209
%
% Fig TEXT object
%
\put{$\xi$} [lB] at  9.620 22.796
\linethickness=0pt
\putrectangle corners at  1.841 23.463 and 11.430 19.241
\endpicture
$$
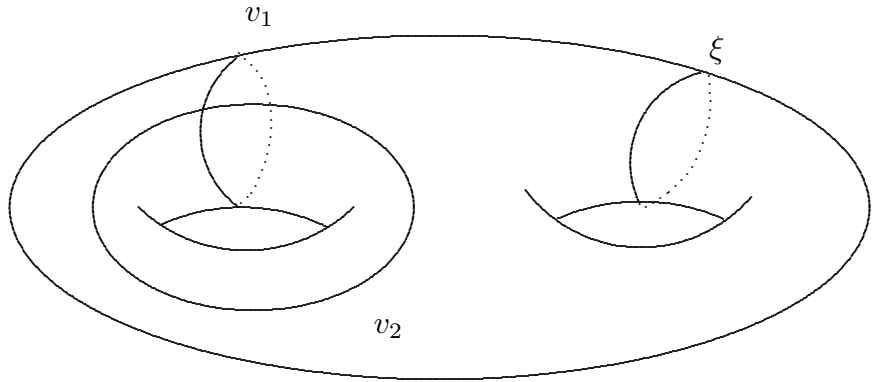
\centerline{Figure 1}
\endinsert
Then by 5.3, $v_1$ and $v_2$ are $\xi$-survivors. Since $v_1\cup v_2\ne 0$, we
obtain $\cl_\k(\xi) =2$. 

Let now $Y$ be an arbitrary closed manifold and $X$ be $\Sigma_g\times Y$. 
Let $\xi'\in H^1(X;\Z)$ be the class with $\xi'|_{\Sigma_g} =\xi$ and $\xi'|_Y=0$. Then clearly
$$\cl_\k(\xi') \ge 2 + \cl_\k(Y),$$
where $\cl_\k(Y)$ is the usual cup-length of $H^\ast(Y;\k)$ (since we may multiply the pullbacks
of $v_1$ and $v_2$ with arbitrary classes coming from $Y$). 

This method produces many examples with large $\cl_k(\xi')$.

\subheading{5.6. Formal spaces} If $X$ is a {\it formal space} \cite{DGMS} 
then all higher Massey products
vanish. We obtain in such a situation that 
{\it a class $v\in H^i(X;\Q)$ is a $\xi$-surviver if and only if $v\cup \xi=0$.} In this case the statement
of Theorem 5.2 is very simple. 

According to \cite{DGMS}, the vanishing of all higher Massey
products takes place in any compact complex manifold $X$
for which the $dd^c$-{\it Lemma} holds (for example, if $X$ is a K\"ahler or a Moi\v sezon space).
This vanishing of higher products directly follows from the diagram
$$\{\E^\ast_M, d\}\overset i\to\leftarrow \{\E^\c_M, d\}\overset \rho\to\to \{H_{d^c}(M),d\}$$
used in the first proof 
of the Main Theorem \cite{DGMS}, cf. page 270.

\subheading{5.7. Remark} {\it Theorem 5.2 becomes false if we allow products (5-4)
with only one $\xi$-surviver instead of two}. Indeed, let $X=S^1\times Y$ and let $\xi\in H^1(X;\Z)$
be the cohomology class of the projection $X\to S^1$. Then $\xi$ can be realized by a closed 
1-form without critical points (the projection). However,
products of the form $\xi\cup v_1\cup\dots \cup v_m$,
where $v_i$ are pullbacks of some classes of $Y$,
which have only one $\xi$-surviver $\xi$, may be nontrivial for $m=\cl_\k(Y)$.

\heading{\bf \S 6. Proof of Theorem 3.1
in the case of integral $\xi$.}\endheading

In this section we will assume that the cohomology class $\xi$ of the form
$\omega$ is integral, $\xi \in H^1(X,\Z)$. We will also assume (without loss of generality) that $\xi$ is indivisible (i.e. not a multiple of an integral class).
Our aim is to prove Theorem 3.1 under these assumptions. The general case will be treated 
later in \S 7.

\subheading{6.1}
There exists a smooth map $f: X \to S^1$ such that $\omega =f^\ast(d\theta)$, where $d\theta$ is the 
standard angular form on the circle $S^1\subset \C$, $S^1=\{z; |z|=1\}.$ 
Suppose that $1\in S^1$ is a regular value of $f$ and denote $f^{-1}(1)$ by $V\subset X$; it is a smooth codimension 
one submanifold. Let $N$ denote the manifold obtained by cutting $X$ along $V$. We will denote by $\partial_+N$
and $\partial_-N$ the components of the boundary of $N$. We obtain a smooth function 
$$g: N \to [0,1],\tag6-1$$
so that 
\roster
\item $g$ is identically 0 on $\partial_+N$ and identically 1 on $\partial_-N$ and has no critical points
near $\partial N$;
\item letting $p: N\to X$ denote the natural identification map, for any $x\in N$,
$$f(p(x)) = \exp(2\pi i g(x));$$
\item $p$ maps the set $S(g)$ of critical points of $g$ homeomorphically onto the set $S(\omega)$
of critical points of $\omega$.
\endroster

\subheading{6.2}
For any subset $X\subset N$ containing $\pn$ 
we will denote by $\cat_N(X|\pn)$ the minimal number $k$ such that $X$
can be covered by $k+1$ closed subsets 
$$X\subset A_0\cup A_1\cup A_2\cup \dots \cup A_k\subset N\tag6-2$$
with the following properties:
\roster
\item $A_0$ is {\it a collar of $\pn$}, i.e. it contains a neighborhood of $\pn$ and 
the inclusion $A_0\to N$ is homotopic to a map $A_0\to \pn$ keeping the points of $\pn\subset A$
fixed;
\item for $j=1, 2, \dots, k,$ each set $A_j$ is disjoint from $\pn$ and the inclusion $A_j\to N$
is null-homotopic.
\endroster

 The number $\cat_N(X|\pn)$ can be viewed as a {\it relative version of the Lusternik
- Schnirelman category}. We refer to \cite{FH, R}, where equivalent (but technically slightly different)
notions were studied.

Standard arguments of the Lusternik - Schnirelman theory give the inequality
$$\cat S(\omega)=\cat S(g) \, \ge \, \cat_N(N|\pn).\tag6-3$$

\subheading{6.3. The deformation complex} In the rest of \S 6 we will establish an inequality
$$\cat_N(N|\pn)\,  \ge \, \cl_{\k}(\xi) - 1.\tag6-4$$
Together with (6-3) this will complete the proof of Theorem 3.1 in the case when $\xi$ is integral.

The proof of (6-4) will consist of building a {\it polynomial deformation} of the cochain complex
$C^\ast(X;a^\xi\otimes F)$ (where $a$ is viewed as a parameter) into 
$C^\ast(N,\pn;F_0)$ "as $a\to \infty$". 
The word {\it "deformation"} is understood here as a finitely generated free cochain complex  $C^\ast$
over the ring $P=\k[\tau]$ of polynomials with coefficients in $\k$ satisfying conditions
(a) and (b) below.
Here $F\to X$ is an arbitrary fixed $\k$-vector bundle over $X$, 
$F_0$ is the induced flat bundle over $N$, and 
$a^\xi$ denotes a flat $\k$-line bundle over $X$ constructed as follows.
Given $a\in \k^\ast$, we will denote by $a^\xi$ the flat $\k$-line bundle over $X$ so that the 
monodromy along any loop $\gamma\in \pi_1(X)$ is $a^{\<\xi, \gamma\>}\in \k^\ast$.

The construction of the deformation $C^\ast$ goes as follows. 
We shall assume that $N$ is triangulated
and $\partial N$ is a subcomplex. Recall that $V=f^{-1}(1)$ (cf. {6.1}) and we will denote 
by $i_\pm : V \to N$ the inclusions, which identify $V$ with $\partial_\pm N$ correspondingly.
The given flat vector bundle $F\to X$ induces the flat vector bundle $F_0\to N$ and also an isomorphism $\sigma: i_+^\ast F_0\to i_-^\ast F_0$.

Denote by $C^q(N)$ and $C^q(V)$ the vector spaces of $F_0$-valued cochains
and $\delta_N: C^q(N) \to C^{q+1}(N)$ and by $\delta_V: C^q(V) \to C^{q+1}(V)$,  
the coboundary homomorphisms. 

Let $C^q(N)[\tau]$ and $C^{q-1}(V)[\tau]$ denote free $P$-modules formed by "polynomials
with coefficients" in the corresponding $\k$-vector spaces; 
for example, an element $c\in  C^q(N)[\tau]$
is a formal sum $c=\sum_{i\ge 0}c_i\tau^i$ with $c_i\in C^q(N)$ and only finitely many $c_i$'s
are nonzero. We shall think of $c\in C^q(N)[\tau]$ as a polynomial "curve",
which associates a point in $C^q(N)$ with any $\tau \in \k$. The $P$-module structure
is given as follows: $\tau\cdot c= \sum_{i\ge 0}c_i\tau^{i+1}$. It is clear that 
$C^q(N)[\tau]$ and $C^{q-1}(V)[\tau]$ are free finitely generated $P$-modules and their ranks
equal to the number of $q$-dimensional simplices in $N$ or $(q-1)$-dimensional simplices
in $V$, correspondingly. 

Consider the natural $P$-module
extensions
$$\delta_N: C^q(N)[\tau] \to C^{q+1}(N)[\tau],\quad\text{ and}\quad 
\delta_V: C^q(V)[\tau] \to C^{q+1}(V)[\tau].\tag6-5$$
They act coefficientwise so that $\delta_N$ and $\delta_V$ are $P$-homomorphisms. 
For example, if $\alpha=\sum_{i\ge 0}\alpha_i\tau^i\in C^q(N)[\tau]$ then 
$\delta_N(\alpha) = \sum_{i\ge 0}\delta_N(\alpha_i)\tau^i$.

Now we define a finitely generated free cochain complex  $C^\ast$
over the ring $P=\k[\tau]$ of polynomials with coefficients in $\k$ as follows:
$C^\ast=\oplus C^q$, where
$$C^q = C^q(N)[\tau]\oplus C^{q-1}(V)[\tau].$$
Elements of $C^q$ will be denoted as pairs $(\alpha,\beta)$, where 
$\alpha\in C^q(N)[\tau]$ and $\beta\in C^{q-1}(V)[\tau]$.
The differential
$\delta: C^q\to C^{q+1}$ is given by the following formula:
$$\delta(\alpha, \beta) = (\delta_N(\alpha), 
(\sigma\otimes i_+^\ast -\tau i_-^\ast)(\alpha) -\delta_V(\beta)),\tag6-6$$
where $\alpha\in C^q(N)[\tau]$ and $\beta\in C^{q-1}(V)[\tau]$. It is clear that $C^\ast$ is the usual
cylinder of the chain map $\sigma\otimes i^\ast_+-\tau i^\ast_-$ with a shifted grading.

We claim now that:{\it 

(a) for any nonzero $a\in \k^\ast$ there is a canonical isomorphism}
$$
\CD
\E^q_a: H^q(C^\ast\otimes _P\k_a) \overset{\simeq}\to\longrightarrow H^q(X;{a^{-\xi}\otimes F}).
\endCD
\tag6-7$$
Here $\k_a$ is $\k$ which is viewed as a $P$-module with the following structure:
$\tau x= ax$ for $x\in \k$.
We will call $\E^q_a$ {\it an isomorphism of evaluation at $\tau=a$};

(b) {\it for $a=0$ we also have a canonical evaluation isomorphism} (although the bundle $a^\xi$
does not exist)
$$\E^q_0: H^q(C^\ast\otimes _P\k_0) \overset{\simeq}\to\longrightarrow H^q(N, \partial_+N;F_0),\tag6-8$$
where $\k_0$ is $\k$ with the following $P$-module structure: $\tau x=0$ for any $x\in \k$.

To show (a) we note that $H^q(X;{a^{-\xi}}\otimes F)$ can be identified with the cohomology of 
complex $C^\ast(X;{a^{-\xi}}\otimes F)$, consisting of cochains $\alpha\in C^q(N)$, satisfying the boundary
conditions
$$a \, i_-^\ast(\alpha) =\sigma\otimes i_+^\ast(\alpha) \in C^q(V).\tag6-9$$
The complex $C^\ast\otimes_P\k_a$ can be viewed as
$$C^q\otimes_P\k_a = C^q(N)\oplus C^{q-1}(V)$$ with the differential given by
$$\delta(\alpha, \beta) = (\delta_N(\alpha), (\sigma\otimes i_+^\ast -a i_-^\ast)(\alpha) -\delta_V(\beta)),\tag6-10$$
where $\alpha\in C^q(N)$ and $\beta\in C^{q-1}(V)$. It is clear that there is a chain homomorphism
$C^\ast(X;{a^{-\xi}}\otimes F)\to C^\ast\otimes_P\k_a$ (acting by $\alpha\mapsto (\alpha,0))$. 
It is easy to see that it induces an isomorphism on the cohomology. Indeed, suppose that a cocycle
$\alpha\in C^q(X;{a^{-\xi}}\otimes F)$ bounds in the complex $C^\ast\otimes_P\k_a$ then there are
$\alpha_1\in C^{q-1}(N)$, $\beta_1\in C^{q-2}(V)$ such that 
$\alpha=\delta_N(\alpha_1)$, $i_+^\ast(\alpha_1) - ai_-^\ast(\alpha_1)-\delta_V(\beta_1)=0.$
We may find a cochain $\beta_2\in C^{q-2}(N)$ such that 
$\sigma\otimes i_+^\ast(\beta_2)=\beta_1$ and 
$i_-^\ast(\beta_2)=0$ (by extending $\beta_1$ into a neighborhood of $\pn$). Then setting
$\alpha_2=\alpha_1-\delta_N(\beta_2)$ we have 
$$\alpha=\delta_N(\alpha_2),\qquad \sigma\otimes i_+^\ast(\alpha_2) - ai_-^\ast(\alpha_2)=0,\tag6-11$$
which means that $\alpha$ also bounds in $C^q(X;{a^{-\xi}}\otimes F)$. 

Similarly, suppose that $(\alpha, \beta)$ is a cocycle of complex $C^\ast\otimes_P\k_a$. As above we
may find a cochain $\beta'\in C^{q-1}(N)$ with $\sigma\otimes i_+^\ast(\beta')=\beta$ and 
$i_-^\ast(\beta')=0$. Then $(\alpha-\delta_N(\beta'), 0)$ is cohomologous to the initial cocycle 
$(\alpha, \beta)$ and is a cocycle of $C^\ast(X;{a^{-\xi}}\otimes F)$.

This proves (a). The statement  (b) follows similarly.\qed

\subheading{6.4. Relative deformation complex} We will define now a relative version of the deformation complex $C^\ast$.

Let $A\subset N$ be a simplicial subcomplex. We will assume that $A$ {\it is disjoint from} $\partial_+N$. Let $C^q(N,A)$ denote the $\k$-vector space of $F_0$-valued cochains on 
$N$, which vanish on $A$. Let $C^q(N, A)[\tau]$ be the polynomial extension
constructed similarly to 
$C^q(N)[\tau]$, cf. subsection 6.3 above. We define the complex $C^\ast_A$ as follows:
$$C^q_A = C^q(N,A)[\tau]\oplus C^{q-1}(V)[\tau].\tag6-12$$
The differential
$\delta: C^q_A\to C^{q+1}_A$ is defined by the following formula:
$$\delta(\alpha, \beta) = 
(\delta_{N,A}(\alpha), (\sigma\otimes i_+^\ast -\tau i_-^\ast)(\alpha) -\delta_V(\beta)),\tag6-13$$
where $\alpha\in C^q(N,A)[\tau]$ and $\beta\in C^{q-1}(V)[\tau]$. 
Here
$\delta_{N,A}: C^q(N,A) \to C^{q+1}(N,A)$ and $\delta_V: C^q(V) \to C^{q+1}(V)$  denote
the coboundary homomorphisms and also the $P$-module extension. 
$i_\pm^\ast : C^q(N,A) \to C^q(V)$
denote the restriction maps of chains, and the same symbols denote also their polynomial extensions
$i_\pm^\ast : C^q(N,A)[\tau] \to C^q(V)[\tau]$. 

Similarly to statements (a) and (b) in 6.3 we have:

{\it (a') for any $a\in \k^\ast$ there is a natural isomorphism
$$H^i(C^\ast_A\otimes_P \k_a) \simeq H^i(X,p(A);{a^{-\xi}}\otimes F),\tag6-14$$
where $p:N \to X$ is the identification map, cf. 6.1;

(b') also,}
$$H^i(C^\ast_A\otimes_P \k_0) \simeq H^i(N, A\cup \pn; F_0).\tag6-15$$
 
\proclaim{6.5. Proposition (The lifting property)} Suppose that $A\subset N$ is a subcomplex 
disjoint from $\pn$ such that the inclusion
$A\to N$ is homotopic to a map with image in a collar of $\pn$. 
Then the homomorphism
$C^\ast_A \to C^\ast$ induces an epimorphism on the cohomology
$$H^i(C^\ast_A\otimes_P \k_a) \to H^i(C^\ast\otimes_P \k_a), \qquad i=0, 1, 2, \dots\tag6-16$$ 
for a generic $a\in \k$ (i.e. for all $a\in \k$, except possibly finitely many points).
\endproclaim
\demo{Proof} Let $\k_0$ denote the field $\k$ considered as a $P$-module with 
trivial $\tau$ action, i.e. $\k_0= P/\tau P$. We will show first that
$$H^i(C^\ast_A\otimes_P \k_0) \to H^i(C^\ast\otimes_P \k_0)\tag6-17$$
is an epimorphism.
We know from  (b') of subsection 6.4 that
$$H^i(C^\ast_A\otimes_P \k_0) \simeq H^i(N, A\cup \pn; F_0)\quad\text{and}\quad 
H^i(C^\ast\otimes_P \k_0) \simeq H^i(N,  \pn;F_0).$$
In the exact sequence
$$\dots \to H^i(N, A\cup \pn;  F_0) \to H^i(N,  \pn;  F_0) 
\overset{j^\ast}\to\longrightarrow 
H^i(A\cup \pn, \pn;F_0)\to \dots $$
$j^\ast$ acts trivially and hence
$H^i(N, A\cup \pn; F_0) \to H^i(N,  \pn; F_0)$ is an epimorphism. This proves that (6-17) 
is an epimorphism.
Now, Proposition 6.5 follows from Lemma 6.6 below. \qed
\enddemo

\proclaim{6.6. Lemma} Let $C$ and $D$ be cochain complexes of free finitely generated 
$P=\k[\tau]$-modules and let $f:C \to D$ be a chain map. Suppose that for some $q$ 
the induced map $f^\ast: H^q(C\otimes_P \k_0) \to H^q(D\otimes_P \k_0)$
is an epimorphism; here $\k_0$ is $\k$ considered with the trivial $P$-action: 
$\k_0= P/\tau P$ . 
Then for a generic $a\in \k^\ast$ the homomorphism
$$f^\ast: H^q(C\otimes_P\k_a) \to H^q(D\otimes_P\k_a)\tag6-18$$
is an epimorphism; here $\k_a$ denotes $\k$ with $\tau$ acting as the multiplication by $a$.\endproclaim 

\demo{Proof} Denote by $Z^q(C), Z^q(D)$ the cocycles of $C$ and $D$ and by $B_q(C)$ and $B^q(D)$ their coboundaries. Recall that $P$ is a PID, and hence
$Z^q(C)$ and $Z^q(D)$ are free $P$-modules. 

Choose bases for $Z^q(C), Z^q(D)$ and $D^{q+1}$ and express in terms of these bases the map
$$f\oplus d: Z^q(C)\oplus D^{q+1} \to Z^q(D).\tag6-19$$
The resulting matrix $M$ is a rectangular matrix with entries in $P$. 

We claim: {\it there exists a minor $A(\tau) \in P$ of the matrix $M$ of size $\rk Z^q(D)\times \rk Z^q(D)$,
such that $A(0)\ne 0$.}
In fact, this claim is clearly {\it equivalent} to the requirement that
$f^\ast: H^q(C\otimes_P \k_0) \to H^q(D\otimes_P \k_0)$ is an isomorphism.

Thus for a generic $a\in \k$ (except finitely many roots of $A(\tau)=0$) it follows that
$A(a)\ne 0$ and hence we obtain that the homomorphism
$$f\oplus d: (Z^q(C)\otimes_P \k_a)\oplus (D^{q+1}\otimes_P \k_a )\to 
Z^q(D)\otimes_P \k_a\tag6-20$$
(which is described by the matrix $M$ with substitution $\tau=a$) 
and hence also (6-18) are epimorphisms. \qed
\enddemo

\proclaim{6.7. Corollary} Let $A\subset X$ be a closed subset such that $A=p(A')$, where $A'\subset N-\pn $ is a closed polyhedral subset such that the inclusion $A'\to N$ is homotopic to a map with values in a collar of
$\pn$. Then the 
restriction map
$$H^q(X,A;{a^\xi}\otimes F)\to H^q(X;{a^\xi}\otimes F)\tag6-21$$
is an epimorphism for almost all $a\in \k^\ast$ (i.e. for all, except finitely many).\endproclaim
\demo{Proof} We just combine the isomorphisms (a) and (a') (cf. 6.3, 6.4) and Proposition 6.5.
\qed \enddemo

\proclaim{6.8. Lemma} Suppose that $E_1$ and $E_2$ are two $\xi$-generic flat bundles over $X$,
and $w\in H_{d_1+d_2}(X;E_1^\ast\otimes E_2^\ast)$ is such that the pairing 
$$\Psi_w: H^{d_1}(X;E_1)\otimes H^{d_2}(X;E_2)\to \k,\tag6-22$$
given by 
$$\Psi_w(v_1\otimes v_2) \, =\, \<v_1\cup v_2,w\>\in \k\tag6-23$$
for $v_i\in H^{d_i}(X;E_i)$, $i=1,2$, is nontrivial (i.e. is not identically zero).
Then for almost all $a\in \k^\ast$ the pairing
$$\Psi^a_w: H^{d_1}(X;a^\xi\otimes E_1)\otimes H^{d_2}(X;a^{-\xi}\otimes E_2)\to \k,\tag6-24$$
acting by formula (6-23) is nontrivial.\endproclaim
The proof will be given in \S 9.6; it will use the preliminary material described in \S 9.1, 9.2.

\subheading{6.9. End of Proof of Theorem 3.1 for integral $\xi$} We need to establish inequality (6-4). 
In other words,
we want to prove the triviality of any cup-product
$$v_1\cup v_2\cup v_3 \cup\dots \cup v_{m+2},\tag6-25$$
where
$v_i\in H^{d_i}(X;E_i)$, $i=1, \dots, m+2$, $d_3>0, \dots, d_{m+2}>0$ and the bundles $E_1$ and $E_2$ are $\xi$-generic. Here
$m$ denotes the relative category $m= \cat_N (N|\pn)$.

Suppose the contrary, i.e. that there exists a nontrivial product (6-25).

We know that $N$ can be covered by closed subsets $A_0\cup A_1\cup  \dots\cup A_m=N$
so that $A_0$ is a collar of $\pn$, and for $j=1,2, \dots, m$ the subset $A_j$ is disjoint from
$\pn$ and null-homotopic in $N$. Hence $p(A_j)$ is null-homotopic in $X$ for $j=1, 2, \dots, m$.
We may assume (without loss of generality)
that the sets $A_j$ are 
polyhedral.  Set $B_j=p(A_{j-2} )$ for $j=3, \dots, m+2$.
We also want to define the sets $B_1$ and $B_2$ as follows.
Let $U_\pm$ be a small cylindrical neighborhood of 
$\partial_\pm N$ in $N$. Set $B_2=p(A_0 -U_+)$.
Let $B_{1}$ be a closed cylindrical neighborhood of $V$ in $X$
containing $\overline {p(U_-)}\cup \overline{ p(U_+)}$. 

Applying Corollary 6.7 combined with Lemma 6.8, we may assume that the two homomorphisms
$$H^{d_1}(X, B_1;E_1) \to H^{d_1}(X;E_1), \quad H^{d_2}(X, B_2;E_2) \to H^{d_2}(X;E_2)
\tag6-26$$
are epimorphisms; if not, we may replace the bundle $E_1$ by $a^\xi\otimes E_1$ and 
the bundle $E_2$ by 
$a^{-\xi}\otimes E_2$ for some $a\in \k^\ast$, leaving the other bundles $E_j$, with $j>2$, 
without modification. Thus we may lift the classes $v_1$ and $v_2$ to classes 
$\tilde v_1\in H^{d_1}(X, B_1;E_1)$ and $\tilde v_2\in H^{d_2}(X, B_2;E_2)$, correspondingly.
On the other hand, 
(since $d_j>0$ for $j=3, \dots, m+2$) we may lift the class
$v_j$ for $j=3, \dots, m+2$ to a relative cohomology class $\tilde v_j\in H^{d_j}(X,B_j;E_j)$.

Now, we get a contradiction to the assumption that the product (6-25) is nontrivial.
Indeed, it is obtained from the product 
$\tilde v_1\cup\tilde v_2\cup \tilde v_3\cup\dots \cup \tilde v_{m+2}$ 
(lying in the space 
$H^d(X, \cup_{j=1}^{m+2} B_j;E_1\otimes E_2\dots \otimes E_m)$, 
where $d=\sum_{j=1}^{m+2} d_j$), by restricting onto $X$, and the cohomology
$H^d(X, \cup_{j=1}^{m+2} B_j;E_1\otimes E_2\dots \otimes E_m)$
vanishes since $X=\cup_{j=1}^{m+2} B_j$. \qed

\heading{\bf \S 7. Proof of Theorem 3.1 in the general case}\endheading

Suppose first that $\xi\in H^1(X;\R)$ {\it has rank 1}, 
i.e. the image of the homomorphism
$H_1(X)\to \R$ determined by $\xi$ is a free abelian group of rank 1 (an infinite cyclic group). Then $\xi$ can be written as
$\xi=\lambda\cdot \xi_0$ where $\xi_0$ is integral, 
$\xi_0\in H^1(X;\Z)$, $\lam\in \R$, and Theorem 3.1 clearly holds for
$\xi$ since (as it was proven) it holds for $\xi_0$. Thus Theorem 3.1 is established for all 
rank 1 classes.

The proof of Theorem 3.1 for classes $\xi$ of rank $>1$ consists of reducing
it to the case of classes of rank 1, which are dense in $H^1(X;\R)$. It is clear that small perturbations
of a closed 1-form may turn it into a rank 1 form with the same set of critical points. 

In order to perform this plan, we need information on the dependence of the cup-length $\cl_\k(\xi)$ on the class 
$\xi\in H^1(X;\R)$. Recall, that 
from our definition of the cup-length in \S 2, it is clear that it depends only on 
$\ker (\xi)\subset H_1(X)$.

Consider $N_\xi = \{\eta\in H^1(X;\R); \eta|_{\ker(\xi)} =0\}$. If $\eta\in N_\xi$ then $\V_\eta\subset
\V_\xi$ and so any $\xi$-generic flat bundle is $\eta$-generic. Hence, $\cl_\k(\eta)\ge \cl_\k(\xi)$.

 Let $\omega$ be a closed 1-form
lying in the cohomology class $\xi\in H^1(X;\R)$ of rank $r>1$. Let $S=S(\omega)$ denote the set
of zeros of $\omega$. 

Choose integral classes $\eta_1, \eta_2, \dots, \eta_r\in N_\xi$ forming a basis of $N_\xi$.
Since $\eta_i|_S =0$, we may realize $\eta_i$ by a smooth closed 1-form $\omega_i$ so that it 
vanishes identically in a neighborhood of $S$; here $i=1, 2, \dots, r$. 
We may write $\xi = \sum_{i=1}^r \alpha_i \eta_i$, with $\alpha_i\in \R$. Let $\tilde \alpha_i$ be a rational number approximating $\alpha_i$, so that $|\alpha_i -\tilde \alpha_i|<\epsilon$. 
Consider the closed 1-form 
$$\tilde \omega = \omega - \sum_{i=1}^r (\alpha_i - \tilde\alpha_i)\omega_i.\tag7-1$$
It coincides with $\omega$ in an open neighborhood $U$ 
of the set $S=S(\omega)$ of critical points of $\omega$. On the other hand,
it has no critical points outside $U$, if $\epsilon >0$ is small enough. Hence, $S(\tilde\omega) =S(\omega)$. The cohomology class 
$\tilde \xi = \sum_{i=1}^r\tilde \alpha_i\eta_i$ of $\tilde \omega$ has rank 1. Hence we obtain
$$\cat(S(\omega)) = \cat (S(\tilde\omega)) \ge \cl_\k (\tilde \xi) -1 \ge \cl_\k (\xi) -1.\tag7-2$$
This completes the proof. \qed

\heading{\bf \S 8. Proof of Theorem 3.2}\endheading

It is enough to prove Theorem 3.2 assuming that the given class $\xi\in H^1(M;\Z)$ is indivisible and $\xi\ne 0$. We may realize the dual homology class to $\xi$
by a connected submanifold $V$, according to a well-known theorem of R. Thom. 
Let us cut $M$ along $V$ and denote 
the resulting cobordism by $N$. 

Using the theorem of Smale \cite{S}, we may find a self-indexing
Morse function $f$ on $N$, so that $f$ assumes constant values on the boundary of 
the cobordism $N$. We may also assume that all the level sets
$f^{-1}(c)$ are connected (otherwise one may use the technique of Smale \cite{S} to perform some cancelations of critical points of indices $0$ and $1$ or critical points of indices $n-1$ and $n$). 

Now, one observes that $f$ has no critical points of indices $0$ and $n$ (since all the level sets 
$f^{-1}(c)$ are connected). Thus, $f$ may have only critical points of indices $1, 2, \dots, n-1$.
Since the level sets are connected, we may collide all the critical points of the same index into a
single (degenerate) critical point, using the method of Takens \cite{T}, \S 2. As a result, we get
a function $g$ on $N$ with the following properties: 
\newline
(a) $g$ has at most $n-1$ critical points in the interior of $N$; \newline
(b) $g$ assumes a constant noncritical value on each boundary component of $N$.
\newline For example, we may assume that $g(N)\subset [-1, n+1]$ and $g^{-1}(-1) = \pn$, 
$g^{-1}(n+1) = \partial_- N$. We obtain the following smooth map into the circle $h: M\to S^1=\R/\Z$,
where $h(x) = \exp(2\pi i g(p(x))/(n+2))$, and $p:N\to M$ is the canonical projection. It is clear that
$h$ has at most $n-1$ critical points. \qed

\heading{\bf \S 9. Proof of Theorem 5.2 and Lemma 6.8}\endheading

Our plan will be as follows. Assuming that there exists a non-trivial product 
$v_1\cup v_2\cup v_3\cup\dots\cup v_m\ne 0$ as in (5-4), so that the classes $v_1$ and $v_2$ are
$\xi$-survivors, we will show that one may deform the classes $v_1$ and $v_2$ to produce
{\it two families of classes $v_1(a)\in H^{d_1}(X;a^{-\xi})$ and 
$v_2(a)\in H^{d_2}(X;a^{\xi})$ "rationally depending on $a\in\k^\ast$ "} so that
for a generic $a\in \k^\ast$ the cup-product
$v_1(a)\cup v_2(a)\cup v_3\cup\dots\cup v_m$ is nontrivial. The condition that $v_1$ and $v_2$
are $\xi$-survivors makes this deformation possible.

We will start with a
general discussion of deformations of chain complexes.

\subheading{9.1. Deformations and their spectral sequences} 
Here we will discuss a spectral sequence
associated with a deformation of a chain complex. The Massey operations (cf. 5.1) appear as
the differentials of a spectral sequence of this type. 

In a much more general situation of deformation of {\it elliptic complexes} the
spectral sequence of deformation was described in \cite{F1}. 

Let $(C, \delta_t)$ be a cochain complex of finitely dimensional $\k$-vector spaces so that
the differential
$\delta_t$ depends polynomially on a parameter $t$. In other words,
$$\delta_t = \delta_0 + t \delta_1 + t^2 \delta_2 + \dots + t^s \delta_s.\tag9-1$$
Here each $\delta_i$ is a degree 1 morphism, $\delta_i\in \Hom(C^\ast, C^{\ast +1})$, and we assume that 
$$\sum_{i+j=k} \delta_i \delta_j =0\qquad\text{for}\quad k=0, 1, \dots.\tag9-2$$
In particular,
$$\delta_0^2 =0,\quad\text{and}\quad\delta_0 \delta_1 + \delta_1 \delta_0 =0.\tag9-3$$
Hence, for any value $t\in \k$, the homomorphism $\delta_t$ is a differential on $C$, 
i.e. $\delta_t^2=0$. 

In the above situation we will describe a spectral sequence $(E_r, d_r)$, where $r=1, 2, \dots, \infty$,
with the following properties:

{\it (i) the initial term is $E_1^\ast =H^\ast(C, \delta_0)$ (cohomology at the point $t=0$).

(ii) the first differential $d_1: E_1^\ast \to E_1^{\ast+1}$ is the homomorphism induced by the
chain map $\delta_1$ (cf. (9-3)).

(iii) the differential $d_r: E_r^\ast \to E_r^{\ast +1}$ depends only on the operators $\delta_0, \delta_1,\dots, \delta_r$ (appearing in (9-1)).

(iv) for all large $r$ the differential $d_r$ vanishes, and 
the limit term $E^\ast_\infty$ is isomorphic to the cohomology $H^\ast (C,\delta_t)$ for a generic $t\in \k$ (i.e. for all $t\in \k$, except finitely many).}

To construct this spectral sequence, consider the cochain complex $(C[t], \delta)$,
where $C[t]= \k[t]\otimes_\k C$ is the space of polynomial curves with values in $C$, and 
the differential $\delta: C[t]\to C[t]$ is a $\k[t]$-homomorphism acting as follows:
$$\delta(c) = \delta_0(c)+t \delta_1(c) +\dots +t^s \delta_s(c),
\quad\text{where}\quad c\in C.\tag9-5$$
We have the following short exact sequence of cochain complexes
$$0\to (C[t], \delta) \overset{t}\to{\to} (C[t], \delta)\to (C,\delta_0)\to 0\tag9-6$$
(multiplication by $t$ and evaluation at $t=0$). The corresponding cohomological long exact sequence
gives an exact couple (cf. \cite{MT})
$$
\CD
H^\ast(C[t],\delta)  \overset{t}\to{\to}  H^\ast(C[t],\delta)\\
\nwarrow \qquad \swarrow\\
 H^\ast(C,\delta_0).
\endCD\tag9-7
$$
The sw-arrow is the evaluation at $t=0$ and the nw-arrow is the boundary homomorphism.
We obtain the spectral sequence generated by this exact couple, cf. chapter 7 in \cite{MT}.
The initial term is $E_1^\ast =H^\ast(C,\delta_0)$ and the first differential $d_1$ clearly coincides with 
the homomorphism induced by $\delta_1$ (using the general definition of the differentials $d_r$, cf. \cite{MT}).

In order to verify (iv), we will first show that
$$\dim_\k E_\infty^i = \rank_{\k[t]} H^i(C[t],\delta).\tag9-8$$
In fact, the $r$-th derived couple 
$$
\CD
D_r \, \, \, \, \overset{t}\to{\to} \, \, \, \, D_r\\
 \nwarrow \, \, \,  \, \, \swarrow\\
E_r
\endCD
$$
has $D_r = t^r H^\ast(C[t],\delta)$. It follows that for large $r$ the multiplication by $t$
(the horizontal arrow) will be monomorphic. Hence, for large $r$ we obtain the short exact sequence
$$0\to t^r H^\ast(C[t],\delta)\overset{t}\to{\to}t^r H^\ast(C[t],\delta)\to E_r^\ast\to 0,
\tag9-9$$
which clearly implies (9-7).

Now we will show that $\rank_{\k[t]} H^i(C[t],\delta)$ coincides with the dimension of the cohomology 
$H^\ast(C, \delta_t)$ for a generic $t\in \k$. For any $a\in \k$ we denote by
$\k_a$ the $\k[t]$-module on the underlying vector space $\k$ with $t:\k \to \k$ acting as 
multiplication by $a\in \k$. The complex $(C, \delta_a)$ (the initial chain complex with the differential
(9-1) where we set $t=a$) coincides with $(C[t],\delta)\otimes_{\k[t]} \k_a$. Hence we may use the
Universal Coefficients Theorem in order to compute the cohomology of $(C, \delta_a)$. We obtain
$$H^i(C,\delta_a) \, =\, H^i(C[t],\delta)\otimes_{\k[t]} \k_a \, \oplus \, 
 \Tor_{\k[t]} (H^{i+1}(C[t],\delta), \k_a).\tag9-10$$
Note that the polynomial ring $\k[t]$ is a PID and so it is clear that the $\Tor$-summand in (9-10)
vanishes for almost all $a\in\k$ (i.e. for all except finitely many). Thus we obtain that for a generic 
$a\in \k$,
$\dim_\k H^i(C, \delta_a) = \rank_{\k[t]} H^i(C[t],\delta)$.

According to the general rules of constructing the spectral sequences out of exact couples 
(\cite{MT}), cf. also \cite{F}, pages 552 - 555, we have 
$$E^i_r = Z_r^i/(t Z^i_{r-1}+t^{1-r}\delta Z^{i-1}_{r-1}),\tag9-11$$
where $Z^i_r$ is defined as the subspace $Z^i_r\subset C[t]$ consisting of chains
$c=c_0+t c_1+ t^2 c_2+\dots$ with 
$$\sum_{i+j=\ell} \delta_i (c_j)=0, \qquad \text{for}\quad \ell =0, 1, \dots, r-1.\tag9-12$$
The differential $d_r: E^i_r \to E^{i+1}_r$ acts as the map $t^{-r}\delta: Z^i_r \to Z^{i+1}_r$.

We may also express the constructed spectral sequence $(E_r, d_r)$ in terms of the ${\k[t]}$-module $H^\ast(C[t], \delta)$, cf.
\cite{F}, pages 552 - 555.
Let $H^\ast(C[t], \delta)=F^i\oplus T^i$ be the representation as the sum of its free and torsion
parts (as modules over ${\k[t]}$). Then
$$E^i_r = F^i/t F^i \oplus T^i/[t T^i + (T^i)_{r-1}] \oplus (T^{i+1})_r/(T^{i+1})_{r-1},
\tag9-13$$
where $(T^i)_{r-1}$ denotes the subspace of $T^i$ consisting of cohomology classes $u$ with 
$t^{r-1}u=0$. 
The action of $d_r$ in terms of this isomorphism is as follows: it vanishes on the first and the
second summand of the decomposition (9-13) and maps the third summand of the decomposition of
$E^i_r$ into the second summand of the decomposition of $E^{i+1}_r$ via the obvious homomorphism
$$(T^{i+1})_r/(T^{i+1})_{r-1}\to (T^{i+1})_r/[t T^{i+1}+(T^{i+1})_{r-1}].\tag9-14$$

From the above computation of the spectral sequence of the deformation 
we obtain the following corollary:

\proclaim{9.2. Corollary} (i) A cohomology class $u\in H^i(C,\delta_0) =E_1^i$ is a cycle of 
all differentials (i.e. for all $r$, $d_r(u)=0$) if and only if there exists
a cohomology class $\tilde u \in H^i(C[t],\delta)$, such that $\tilde u|_{t=0} =u$.
\newline
(ii) A cohomology class $u\in H^i(C,\delta_0) =E_1^i$ survives up to 
$E_\infty^i$ (i.e. for all $r$, $d_r(u)=0$ and $u\notin \im(d_r)$) if and only if there exists
a cohomology class $\tilde u \in H^i(C[t],\delta)$, such that $\tilde u|_{t=0} =u$
and there are no classes $\tilde u\in H^i(C[t],\delta)$, such that 
$\tilde u|_{t=0} =u$ and $t^r\tilde u=0$ for some $r\ge 1$.\newline
(iii) the following conditions are equivalent: \newline

(a) all the differentials $d_r$ are trivial, $d_r\equiv 0$;\newline

(b) the $\Cal R$-modules $H^i(C[t],\delta)\otimes_{\k[t]}{\Cal R}$ are free, where 
${\Cal R}$ denotes the ring of rational functions $p(t)/q(t)$, where $q(0)\ne 0$, and $p(t), q(t)\in
\k[t]$. 

(c) for any fixed $i$, the dimension of the cohomology $H^i(C,\delta_0)$ equals the minimum
of the dimensions
$\dim_\k H^i(C,\delta_a)$, where $a\in \k$. 
\endproclaim

Indeed, in the notation above, we have 
$$E^i_1 = F^i/t F^i \oplus T^i/t T^i  \oplus (T^{i+1})_1.\tag9-15$$
A class $u\in E_1^i$ is a cycle of all differentials iff it has a trivial $(T^{i+1})_1$-component.
A class $u\in E_1^i$
survives up to $E_\infty^i$ if and only if it has a nontrivial $F^i/t F^i $
component. 

In this paper we will actually deal with a special kind of deformation (9-1) having the form $\delta_t = \delta_0 +t\delta_1$
({\it linear deformation}). In this case some general formulae become simpler. For example, 
here $Z^i_r$ is the subspace $Z^i_r\subset C[t]$ consisting of chains
$c=c_0+t c_1+ t^2 c_2+\dots$ with 
$$\delta_0(c_\ell) +\delta_1(c_{\ell -1})=0,\qquad \text{for}\quad \ell =0, 1, \dots, r-1.\tag9-17$$
The differential $d_r: E^i_r \to E^{i+1}_r$ acts as follows
$$d_r(c)\, \, =\, \,  \delta_1(c_{r-1}).\tag9-18$$

\subheading{Remark} Given a deformation $(C,\delta_t)$, one may replace the parameter $t$ by
$\tau=t-a$. In general, the obtained spectral sequence will be different, but will have the same limit 
term. We will call it {\it the spectral sequence of the deformation centered at} $t=a$. 

\subheading{9.3} Now we will apply the general facts concerning deformations of chain complexes
(described above in 9.1 - 9.2),
to the special situation of {\it the deformation complex} 6.3. 

We will adopt the notation introduced in 6.3. 
In particular, $X$ denotes a closed manifold and $V=f^{-1}(1)$ (cf. {6.1}) is a codimension
one submanifold realizing the class $\xi$; $N$ is obtained by cutting $X$ along $V$.
We assume that $N$ is triangulated
and $\partial N$ is a subcomplex. The maps $i_\pm : V \to N$ denote the inclusions.
$C^q(N)$ and $C^q(V)$ denote the vector spaces of $\k$-valued cochains
and $\delta_N: C^q(N) \to C^{q+1}(N)$ and $\delta_V: C^q(V) \to C^{q+1}(V)$ denote
the coboundary homomorphisms. 

Now we define {\it the deformation cochain complex} $(C^\ast, \delta_t)$ as follows:
$C^\ast=\oplus C^q$, where
$C^q = C^q(N)\oplus C^{q-1}(V).$
Elements of $C^q$ will be denoted as pairs $(\alpha,\beta)$, where 
$\alpha\in C^q(N)$ and $\beta\in C^{q-1}(V)$.
The differential
$\delta_t: C^q\to C^{q+1}$ is given by the following formula:
$$\delta_t(\alpha, \beta) = (\delta_N(\alpha), 
( i_+^\ast -(1+t) i_-^\ast)(\alpha) -\delta_V(\beta)),\tag9-19$$
where $\alpha\in C^q(N)$ and $\beta\in C^{q-1}(V)$. 

Comparing with formula (9-5) we see that this deformation is linear: 
$\delta_t=\delta_0+t\delta_1$. 

The differential
$\delta_0$ computes the cohomology $H^\ast(X;\k)$ (cf. (6-7)). The deformation
$\delta_1: C^\ast\to C^{\ast+1}$ acts as follows: $\delta_1(\alpha, \beta)=(0, -i_-^\ast(\alpha))$.
We will denote the homomorphism $\delta_1$ by 
$(\alpha,\beta)\mapsto (-1)^{|\alpha|+1}(\alpha,\beta)\cup \xi$, 
where $(\alpha, \beta)\in C^\ast$. 
We will see below in section 9.4 (cf. (9-25)) that this notation is justified.

Applying the spectral sequence of 9.1 to the present situation and using isomorphisms (6-7), (6-8), 
we obtain that: \newline
{\it there exists a spectral sequence $E_r^\ast$, $r\ge 1$, with the initial term $E_1^i=H^i(X;\k)$, such that $\dim_\k E_\infty^i$ coincides with 
$\dim_\k H^i(X;a^\xi)$ for a generic $a\in k^\ast$.}

According to 9.3, we see that the classes in $E_r^i$ of this spectral sequence can be described
as cycles 
$c_0\in C^i$, $\delta_0(c_0)=0$ (considered up to a certain equivalence relation), 
such that one may find cochains $c_1, \dots, c_{r-1}\in 
C^i$ with
$$\delta_0(c_1) = (-1)^i  c_0\cup\xi, \, \, \delta_0(c_2) = (-1)^i  c_1\cup\xi, \dots,\, \, 
\delta_0(c_{r-1}) = (-1)^i c_{r-2}\cup\xi.\tag9-20$$
The differential $d_r:E_r^i\to E_r^{i+1}$ of the spectral sequence  (cf. 9.3)
maps the class of $c_0$ into
the class of $(-1)^{i+1}c_{r-1}\cup\xi$. Comparing this description with \cite{N3}
we observe that {\it the differentials of this spectral sequence are precisely
the higher Massey operations described by S. P. Novikov} \cite{N3}. Pazhitnov \cite{P} showed that
these Massey operations are equivalent to the {\it symmetric Massey products}, 
which were introduced
and studied by D. Kraines \cite{K}.

\subheading{9.4. Multiplication in the deformation complex} We will now describe the multiplication
$$\psi_t: C^\ast\otimes_\k C^\ast \to C^\ast,\tag9-21$$
where $C^\ast$ is the deformation complex with the differential $\delta_t$ 
as in 9.3. Note that $\psi_t$ depends linearly on the parameter $t$.

Given $c=(\alpha, \beta)\in C^q$ and
$c'=(\alpha', \beta')\in C^{q'},$ we set
$$\psi_t(c \otimes c') = (\alpha\cup \alpha', (-1)^{|\alpha|}(1+t)i_-^\ast(\alpha)\cup\beta'+\beta\cup i_+^\ast(\alpha'))\in C^{q+q'}.\tag9-22$$
It is not difficult to check (an unpleasant exercise!) that 
$$
\delta_0[\psi_t((\alpha, \beta)\otimes (\alpha',\beta'))] = \psi_t(\delta_t(\alpha,\beta)\otimes (\alpha',\beta')) + (-1)^{|\alpha|}
\psi_t((\alpha,\beta)\otimes \delta_{t'}(\alpha',\beta')),\tag9-23$$
where $t' = (1+t)^{-1}-1$. Hence, the pairing (9-21), viewed as a chain map
$$\psi_t: (C^\ast, \delta_t)\otimes_\k (C^\ast, \delta_{t'})\to (C^\ast, \delta_0),$$
satisfies the Leibnitz rule.

If $a\in \k^\ast$ and $t=a-1$ then by (6-7) we know that $H^i(C^\ast,\delta_t)\simeq H^i(X;a^{-\xi})$ and
$H^j(C^\ast,\delta_{t'})\simeq H^j(X;a^{\xi})$, where $t' = (1+t)^{-1}-1$. Also, 
$H^\ast(C^\ast,\delta_{0})\simeq H^\ast(X;\k)$, again by (6-7). It is quite obvious that {\it the
cohomological product 
$$H^i(C^\ast,\delta_t)\otimes H^j(C^\ast,\delta_{t'})\to H^{i+j}(C^\ast,\delta_{0}),$$
derived from the pairing (9-21), 
coincides via the above mentioned isomorphisms with the cup-product}
$$\cup:\, \, H^i(X;a^{-\xi})\otimes H^j(X;a^{\xi})\to H^{i+j}(X;\k).\tag9-24$$
To see this one repeats the arguments used at the end of section 6.3. We may realize each cycle
$(\alpha, \beta)$ of $(C^\ast, \delta_a)$ by a cycle of the form $(\alpha, 0)$, and then
viewing formula
(9-22) makes our statement obvious.

It is clear that the pair $(0,1)\in C^1$ (where $1\in C^0(V;\k)$) is a cycle,
representing the class $\xi\in H^1(X;\k)$ (more precisely, its 
image under the coefficient homomorphism $\Z\to \k$, mapping 1 to 1). Hence we obtain (comparing
with (9-22))
$$\psi_0((\alpha,\beta)\otimes (0,1)) = (0, (-1)^{|\alpha|}i_-^\ast(\alpha)) =  (-1)^{|\alpha|+1}\delta_1(\alpha, \beta),$$
(where $(\alpha,\beta)\in C^\ast$), 
which may be rewritten as
$$\delta_1(c) = (-1)^{i+1} c\cup\xi,\qquad c\in C^i \tag9-25$$
in accordance with our notation used in 9.3.

\subheading{9.5. Completing the proof of Theorem 5.2} 
Assume that we are in conditions of Theorem 5.2, i.e. there exists a non-trivial product 
$v_1\cup v_2\cup v_3\cup\dots\cup v_m\ne 0$ as in (5-4) with $v_1\in H^{d_1}(X;\k)$ and 
$v_2\in H^{d_2}(X;\k)$ being $\xi$-survivors and for $i=3, \dots, m$ the class
$v_i\in H^{d_i}(X;E_i)$ has positive dimension $d_i\ge 1$, where $E_i$ is a flat bundle over $X$.
We would like to show that {\it for a generic} $a\in \k^\ast$ there is a non-trivial cup-product
$v'_1\cup v'_2\cup v_3\cup\dots\cup v_m\ne 0$, where $v'_1\in H^{d_1}(X;a^\xi)$ and
$v'_2\in H^{d_2}(X;a^{-\xi})$.

Non-triviality of the product $v_1\cup v_2\cup v_3\cup\dots\cup v_m$ means that 
there exists a homology class $z\in H_d(X; E^\ast)$, 
where $d=\sum d_i$ and $E=E_3\otimes \dots E_m$, so that 
$\<v_1\cup v_2\cup v_3\cup\dots\cup v_m,z\>\ne 0$. Hence, setting 
$z'=(v_3\cup\dots\cup v_m)\cap z\in H_{d_1+d_2}(X;\k),$ we will have
$\<v_1\cup v_2, z'\>\ne 0$. Realizing $z'$ by a cycle in the complex $(C^\ast, \delta_0)$
(which we will still denote by $z'$), and composing with the multiplication map $\psi_t$ gives 
a family of chain maps
$$\Psi_t: (C^\ast, \delta_t)\otimes _\k (C^\ast, \delta_{t'})\to (\k,0), \qquad \Psi_t(c\otimes c')=
\<\psi_t(c\otimes c'),z'\>,\tag9-26$$
where $(\k,0)$ is a cochain complex with trivial differential having $\k$ in dimension $d_1+d_2$
and all other chain groups are trivial. 
This chain map $\Psi_t$ depends polynomially on the parameters $t$ and $t'=(1+t)^{-1}-1$ and so
can be viewed as assuming its values in the following ring of rational functions 
$(1+t)^{-1}\k[t]$.

Since $v_1$ and $v_2$ are $\xi$-survivors, by Corollary 9.2 we may find polynomial cycles
$$c_t=c_0+tc_1+t^2c_2+\dots, \qquad c'_{t'}=c'_0+t'c'_1+{t'}^2c'_2+\dots$$
with $c_j\in C^{d_1}$ and $c'_j\in C^{d_2}$,
such that $c_0$ represents $v_1$, $c'_0$ represents $v_2$ and the following holds:
$$\delta_t(c_t)=0,\qquad \delta_{t'}(c'_{t'})=0.$$
We obtain that $\Psi_t(c_t\otimes c'_{t'})\in (1+t)^{-1}\k[t]$ is a rational function of $t$, regular
for all $t\ne -1$. Our assumptions give that $\Psi_t(c_t\otimes c'_{t'})\ne 0$ for $t=0$
(since $\<v_1\cup v_2, z'\>\ne 0$). Hence
for almost all $a\in \k^\ast$ the function $\Psi_t(c_t\otimes c'_{t'})$ is nonzero for $t=a-1$. Now, 
evaluating the cycle $c_t$
at $t=a-1$ we obtain a cohomology class $v_1(a)\in H^{d_1}(X;a^{-\xi})$ and 
evaluating the cycle $c'_{t'}$
at $t'=a^{-1}-1$ we obtain a cohomology class $v_2(a)\in H^{d_2}(X;a^{\xi})$ (in view of (6-7)).
We have
$$\<v_1(a)\cup v_2(a),z'\> = \<v_1(a)\cup v_2(a)\cup v_3\cup\dots \cup v_m, z\>\ne 0,$$
for almost all $a\in\k^\ast$ since this value can be obtained by evaluating the rational function $\Psi_t(c_t\otimes c'_{t'})$
at $t=a$.
This completes the proof. \qed

\subheading{9.6. Proof of Lemma 6.8} We will use arguments based on Corollary 9.2,
which are similar to those used in the Proof  of Theorem 5.2 given above.

Consider a local system $\tau^\xi$ over $X$, which has fiber 
$\Lambda = \k[\tau, \tau^{-1}]$ and where
the monodromy along any loop $\gamma\in \pi_1(X)$ is given by multiplication on 
$\tau^{\<\xi,\gamma\>}: \Lambda \to \Lambda$. For $i=1, 2$, the cohomology 
$H^\ast(X; \tau^\xi\otimes_\k E_i)$ is a finitely generated $\Lambda$-module. Its torsion part is supported on a 
finite set of points. Since we assume, that the flat bundles $E_1$ and $E_2$ are $\xi$-generic,
from statement (iii) of Corollary 9.2 we obtain that the point
$\tau=1$ does not belong to this support.

The $\tau=a$ evaluation map $\Lambda\to \k$, where $a\in\k^\ast$, 
$p\in \Lambda\mapsto p(a)$, makes $\k$ a $\Lambda$-module, which 
we will denote by $\k_a$. Using the
Universal Coefficients Theorem we obtain 
$$H^{d_i}(X; E_i)\, \simeq\, H^{d_i}(X;\tau^\xi\otimes_\k E_i)\otimes_\Lambda \k_1\qquad i=1,2,\tag9-27$$ 
(evaluation at $\tau=1$),
where we have used the assumption that $E_i$ are $\xi$-generic.
Hence we conclude that for any pair of classes $v_i\in H^{d_i}(X; E_i)$, $i=1, 2$, 
there exist classes $\tilde v_1\in H^{d_1}(X;\tau^\xi\otimes_\k E_1)$ and 
$\tilde v_2\in H^{d_2}(X;\tau^{-\xi}\otimes_\k E_2)$, such that their evaluation at $\tau=1$
gives $v_i$, i.e. $\tilde v_i(1) =v_i$.

For any $a\in \k^\ast$, we may evaluate $\tilde v_i$ at $\tau=a$ obtaining
$\tilde v_i(a) \in H^{d_i}(X;a^\xi\otimes_\k E_i)$. 

We have the following family of pairings of local systems
$$\Phi^a: (\tau^\xi\otimes_\k E_1)\otimes_\k (\tau^{-\xi}\otimes_\k E_2)\to E_1\otimes E_2, \quad a\in \k^\ast,\tag9-28$$
which act by $(p\otimes e_1)\otimes (q\otimes e_2) \mapsto p(a)q(a)e_1\otimes e_2$, where 
$p,q\in \Lambda$, $e_1\in E_1$, $e_2\in E_2$. Using $\Phi^a$, one constructs a family of cup-products
$$\cup^a: H^{d_1}(X;\tau^\xi\otimes_\k E_1)\otimes H^{d_2}(X;\tau^{-\xi}\otimes_\k E_2)\to H^{d_1+d_2}(X;E_1\otimes E_2).\tag9-29$$
It is clear that for fixed $\tilde v_1\in H^{d_1}(X;\tau^\xi\otimes_\k E_1)$ and 
$\tilde v_2\in H^{d_2}(X;\tau^{-\xi}\otimes_\k E_2)$ the product 
$\tilde v_1\cup^a \tilde v_2\in H^{d_1+d_2}(X;E_1\otimes E_2)$ (as a function of $a$) is a polynomial in $a$ and $a^{-1}$.
We obtain for $a\in\k^\ast$
$$\<\tilde v_1\cup^a \tilde v_2, w\>  = \Psi_w^a(\tilde v_1(a)\otimes \tilde v_1(a)),\tag9-30$$
and so 
for a pair of classes $v_i\in H^{d_i}(X;E_i), i=1, 2$, as above, the function 
$a\mapsto \Psi^a_w( \tilde v_1(a)\otimes \tilde v_2(a))$, where $a\notin S$, 
is a Laurent polynomial in $a$.
This polynomial is nontrivial for $a=1$; hence for almost all $a\in \k^\ast$ we have
$\Psi^a_w( \tilde v_1(a)\otimes \tilde v_2(a))\ne 0$. This completes the proof. \qed

\enddocument